\documentclass{proc-l}

\usepackage{amsmath,amsthm, latexsym, amssymb,mathrsfs,  tikz-cd}
\usepackage{stackengine,scalerel, MnSymbol}
\usepackage[colorlinks]{hyperref}
\usepackage[capitalize]{cleveref}
\definecolor{dark-red}{rgb}{0.6,0,0}
\definecolor{dark-green}{rgb}{0,0.4,0}
\definecolor{medium-blue}{rgb}{0,0,0.5}
\hypersetup{
    colorlinks, linkcolor={dark-red},
    citecolor={dark-green}, urlcolor={medium-blue}
}

\NewDocumentCommand{\liminv}{o}{\IfNoValueTF{#1}{\displaystyle\lim_{ \longleftarrow }}{\displaystyle\lim_{\begin{subarray}{c} \longleftarrow \\[-2pt]  #1 \end{subarray}}}}

\usepackage{todonotes}
\todostyle{commentALL}{
    size=\footnotesize
}
\todostyle{commentMC}{
    commentALL,
    backgroundcolor=gray!5,
    textcolor=blue,
    bordercolor=gray!10,
}
\todostyle{commentSH}{
    commentALL,
    backgroundcolor=brown!10,
    textcolor=red,
    bordercolor=brown!15,
}

\newcommand{\Fil}{\mr{Fil}}

\newcommand{\gr}{\mr{gr}}
\newcommand{\Gr}{\mathrm{Gr}}

\newcommand{\Fl}{\mr{Fl}}

\newcommand{\Spd}{\mathrm{Spd}}
\newcommand{\GL}{\mathrm{GL}}

\newcommand{\mc}[1]{\mathcal{#1}}
\newcommand{\mbb}[1]{\mathbb{#1}}
\newcommand{\mr}[1]{\mathrm{#1}}

\newcommand{\et}{\mathrm{\acute{e}t}}

\newcommand{\dR}{\mathrm{dR}}

\newcommand{\Spa}{\mathrm{Spa}}

\DeclareMathOperator{\Lie}{Lie}
\DeclareMathOperator{\Spec}{Spec}

\DeclareMathOperator{\Gal}{Gal}

\newcommand{\SL}{\mathrm{SL}}
\numberwithin{equation}{subsection}
\numberwithin{equation}{subsubsection}

\theoremstyle{plain}

\newtheorem*{theorem*}{Theorem}

\newtheorem{theorem}[subsubsection]{Theorem}

\newtheorem{proposition}[subsubsection]{Proposition}
\newtheorem{lemma}[subsubsection]{Lemma}

\theoremstyle{definition}

\newtheorem{example}[subsubsection]{Example}

\newtheorem{definition}[subsubsection]{Definition}
\newtheorem{remark}[subsubsection]{Remark}
\newtheorem{question}[subsubsection]{Question}

\newcommand{\Res}{\mathrm{Res}}

\title{Transitivity of the $\mathbb{B}^+_\dR$-loop group action on Schubert cells}
\author{Sean Howe}
\address{University of Utah Department of Mathematics, 155 S 1400 E, Salt Lake City, UT 84112}
\email{sean.howe@utah.edu}

\begin{document}

\begin{abstract}
For $G$ a connected linear algebraic group over a $p$-adic field, we show that the action of $G(\mathbb{B}^+_\dR)$ on each Schubert cell in the $\mathbb{B}_\dR^+$-affine Grassmannian is transitive in the \'{e}tale topology on affinoid perfectoids, generalizing a result in the reductive case due to Fargues and Scholze. 
\end{abstract}

\maketitle

\newcommand{\AffPerfd}{\mathrm{AffPerfd}}
\section{Introduction}
 For $L$ a field and $G/L$ a connected linear algebraic group, the affine Grassmannian of $G$ is the functor on $L$-algebras defined as the \'{e}tale sheafification of $R \mapsto G\left(R((t))\right)/G(R[[t]])$.
It plays a fundamental role in geometric representation theory, e.g. via the geometric Satake isomorphism. 

When $L$ is a $p$-adic field (a complete discretely valued extension of $\mathbb{Q}_p$ with perfect residue field), there is a variation on this functor in $p$-adic geometry obtained by substituting a twisted form of the construction $R \mapsto R[[t]] \subseteq R((t))$: For any perfectoid ring $R$ over $L$, we have Fontaine's canonical deformation of $R$ over $L$, $\theta:\mathbb{B}^+_\dR(R)\twoheadrightarrow R$. The kernel of $\theta$ is a principal ideal and, for any generator $\xi$,  $\mathbb{B}^+_\dR(R)$ is $\xi$-adically complete and $\xi$-torsion free and we write $\mathbb{B}_\dR(R)=\mathbb{B}^+_\dR(R)[1/\xi]$ (these constructions will be recalled in more detail in \cref{s.preliminaries}). The assignment $R \mapsto \mathbb{B}^+_\dR(R) \subseteq \mathbb{B}_\dR(R)$ behaves in many ways like $R \mapsto R[[t]] \subseteq R((t))$ --- indeed, for $K/L$ a perfectoid field, $\mathbb{B}^+_\dR(K)$ is even  isomorphic to $K[[\xi]]$ (as an $L$-algebra but non-canonically). Thus, for $\AffPerfd_L$ the category of affinoid perfectoid spaces over $L$ and $G/L$ a connected linear algebraic group, it is natural to define the $\mathbb{B}^+_\dR$-affine Grassmannian $\Gr_G$ on $\AffPerfd_L$ as the \'{e}tale sheafification\footnote{For $G$ reductive, \cite[Theorem 3.1]{CesnaviciusYoucis.TheAnalyticTopologySufficesForTheBplusDrGrassmannian} shows the sheafification in the analytic topology agrees; in \cref{prop.aff-grass-properties} we explain how to extend this result to non-reductive $G$.} of 
\[ \Spa(R,R^+)/\Spa L \mapsto G(\mathbb{B}_\dR(R)) / G(\mathbb{B}^+_\dR(R)).\]
It plays an important role in the geometrization of the local Langlands correspondence \cite{FarguesScholze.GeometrizationOfTheLocalLanglandsCorrespondence} and as a period domain in $p$-adic Hodge theory  \cite{Scholze.pAdicHodgeTheoryForRigidAnalyticVarieties,ScholzeWeinstein.BerkeleyLecturesOnPAdicGeometryAMS207, HoweKlevdal.AdmissiblePairsAndpAdicHodgeStructuresIITheBiAnalyticAxLindemannTheorem}.  

For $[\mu]$ a conjugacy class of geometric cocharacters of $G$ with reflex field $L([\mu])$ (see \cref{ss.schubert-cells}), the affine Schubert cell $\Gr_{[\mu]} \subseteq \Gr_G \times_{\Spd L} \Spd L([\mu])$ consists of those sections whose pullback to any geometric point $\Spa(C,C^+)$ lies in
\[ G(\mathbb{B}^+_\dR(C)) \xi^{\mu} G(\mathbb{B}^+_\dR(C)) / G(\mathbb{B}^+_\dR(C)) \subseteq G(\mathbb{B}_\dR(C))/G(\mathbb{B}^+_\dR(C)) =\Gr_G( (C,C^+))\]
for one (equivalently, any), choice of $\mu \in [\mu]$ and generator $\xi$ of $\mathrm{Ker}\,\theta$. 
 We refer to \cite[\S5.2]{HoweKlevdal.AdmissiblePairsAndpAdicHodgeStructuresIITheBiAnalyticAxLindemannTheorem} for the basic properties of these affine Schubert cells in the generality of an arbitrary connected linear algebraic group $G/L$.

\begin{remark}\label{remark.local-shimura}In the case of $G$ reductive, these Schubert cells geometrize the Cartan decomposition on geometric points and play a key role, e.g., in the geometric Satake of \cite[Chapter VI]{FarguesScholze.GeometrizationOfTheLocalLanglandsCorrespondence}. In general, Schubert cells do not exhaust the geometric points of $\Gr_G$, but the geometric points lying in a Schubert cell still play a distinguished role: by \cite[Theorem D]{HoweKlevdal.AdmissiblePairsAndpAdicHodgeStructuresITranscendenceOfTheDeRhamLattice}, they classify those lattices such that the associated Bialynicki-Birula filtration is exact. Moreover, as explained in \cite{HoweKlevdal.AdmissiblePairsAndpAdicHodgeStructuresITranscendenceOfTheDeRhamLattice, HoweKlevdal.AdmissiblePairsAndpAdicHodgeStructuresIITheBiAnalyticAxLindemannTheorem}, even if one starts with a reductive $G$, the Schubert cells in $\Gr_H$ for non-reductive subgroups $H \leq G$ arise naturally when describing special subvarieties of the local Shimura varieties (and their non-minuscule analogs) that are associated to $G$.
\end{remark}

By definition, the natural left multiplication action of $G(\mathbb{B}_\dR)$ on $\Gr_{G}$ is transitive in the \'{e}tale topology. The action of $G(\mathbb{B}^+_\dR) \subseteq G(\mathbb{B}_\dR)$ preserves $\Gr_{[\mu]} \subseteq \Gr_G$ and, for $(C,C^+)$ an algebraically closed perfectoid field over $L$, $G(\mathbb{B}^+_\dR(C))$ acts transitively on $\Gr_{[\mu]}(\Spa(C,C^+))$. It is not clear from this definition, however, whether or not $G(\mathbb{B}^+_\dR)$ acts transitively on $\Gr_{[\mu]}$ in the \'{e}tale topology (or even in the $v$-topology), and this property was left as an open question in \cite{HoweKlevdal.AdmissiblePairsAndpAdicHodgeStructuresIITheBiAnalyticAxLindemannTheorem} (see Question 3.2.2 in v1 on ar{X}iv). Our main result is that this action is indeed transitive. \emph{In other words, $\Gr_{[\mu]}$ is an \'{e}tale homogeneous space for the action of  $G(\mathbb{B}^+_\dR)$.}

\begin{theorem}\label{theorem.action-transitive}
    Let $L$ be a $p$-adic field with algebraic closure $\overline{L}$, let $G/L$ be a connected linear algebraic group, and let $[\mu]$ be a conjugacy class of cocharacters of $G_{\overline{L}}$. The action of $G(\mathbb{B}^+_\dR)$ on $\Gr_{[\mu]}$ is transitive in the \'{e}tale topology on $\AffPerfd_{L([\mu])}$: that is, if $s_i \in \Gr_{[\mu]}(\Spa(R,R^+))$ for $i=1,2$, then there is an \'{e}tale cover $\Spa(S,S^+) \rightarrow \Spa(R,R^+)$ and $g \in G(\mathbb{B}^+_\dR(S))$ such that 
    \[ g \cdot s_1|_{\Spa(S,S^+)} = s_2|_{\Spa(S,S^+)}. \]
\end{theorem}

In the case of $G$ reductive, the analog of \cref{theorem.action-transitive} for the $v$-topology is contained in \cite[Proposition VI.2.4]{FarguesScholze.GeometrizationOfTheLocalLanglandsCorrespondence}\footnote{This statement seems to be well-known to experts in the geometrization of the local Langlands correspondence but less well-known to those who are primarily interested in the $\mathbb{B}^+_\dR$-affine Grassmannian as a period domain in $p$-adic Hodge theory, for whom the primary reference is often the earlier work \cite{ScholzeWeinstein.BerkeleyLecturesOnPAdicGeometryAMS207} in which this statement is not established. In particular, one purpose of this note is simply to popularize the statement in the reductive case in a simpler notational setting.}.  To prove \cref{theorem.action-transitive}, we first adapt and expand on\footnote{In particular, \cref{prop.closed-schubert-functor} fills in a point that was not explained in \cite[Proof of Proposition VI.2.4]{FarguesScholze.GeometrizationOfTheLocalLanglandsCorrespondence} and \cref{prop.v-quotient-is-etale-quotient} refines the result from the $v$-topology to the \'{etale} topology.}
the argument in \cite[Proof of Proposition VI.2.4]{FarguesScholze.GeometrizationOfTheLocalLanglandsCorrespondence} to obtain the result for $G$ reductive (\cref{prop.red-case}). We then deduce the general case using arguments similar to those that were used in \cite[\S5.2]{HoweKlevdal.AdmissiblePairsAndpAdicHodgeStructuresIITheBiAnalyticAxLindemannTheorem} to deduce other properties of $\Gr_{[\mu]}$ in the general case from those established in the reductive case in \cite{ScholzeWeinstein.BerkeleyLecturesOnPAdicGeometryAMS207, FarguesScholze.GeometrizationOfTheLocalLanglandsCorrespondence}.

The reductive case of \cref{theorem.action-transitive} is the most important for the applications to the geometrization of the local Langlands correspondence in \cite{FarguesScholze.GeometrizationOfTheLocalLanglandsCorrespondence}, but the non-reductive case arises naturally in the study of special subvarieties of local Shimura varieties (see \cref{remark.local-shimura}) and has applications to the study of a differential structure on Schubert cells even in the reductive case (see \cref{ss.application}).

\begin{example}\label{remark.homogeneous-space}Suppose that $[\mu]$ admits\footnote{This holds, e.g., if $G$ is quasi-split over $L([\mu])$: indeed, write $U$ for the unipotent radical of $G$, and fix a Levi decomposition $G=MU$ defined over $L$ so that $M=G/U$. We claim that $G$ and $M$ have the same conjugacy classes of cocharacters over any extension of $L$: indeed, any cocharacter of $G$ is conjugate to one factoring through $M$, e.g. by conjugating into a maximal torus lying in $M$, and, by projecting to $M$ as a quotient, we see that two cocharacters of $M$ conjugate in $G$ are also conjugate in $M$. In particular, we obtain the same reflex field if we view $[\mu]$ as a conjugacy class of cocharacters of $M$, thus the existence of a representative in the general quasi-split case follows from the existence in the reductive quasi-split case, which is \cite[Lemma 1.1.3-(a)]{Kottwitz.ShimuraVarietiesAndTwistedOrbitalIntegrals}.} a representative $\mu$ defined over $L([\mu])$.
In this case $L([\mu])=L$ and there is a natural point $\ast_{\mu} \in \Gr_{[\mu]}(L)$
(see \cref{ss.canonical-point-and-stabilizer}).  \cref{theorem.action-transitive} is then equivalent to the statement that the orbit map for $\ast_\mu$ induces an isomorphism of \'{e}tale sheaves on $\AffPerfd_{L}$, 
\[ G(\mathbb{B}^+_\dR) / \mathrm{Stab}_{G(\mathbb{B}^+_\dR)}(\ast_\mu)  = \Gr_{[\mu]}, \]
where on the left the quotient is of sheaves for the \'{e}tale topology. We note that, for $\Spa(R,R^+) \in \AffPerfd_L$ and $\xi \in \mathbb{B}^+_\dR(R)$ a generator of $\mr{Ker}\,\theta$,
\[ \mathrm{Stab}_{G(\mathbb{B}^+_\dR(R))}(\ast_\mu) = \xi^{\mu}G(\mathbb{B}^+_\dR(R))\xi^{-\mu}\cap G(\mathbb{B}^+_\dR(R)), \]
where the intersection is taken in $G(\mathbb{B}_\dR(R))$. The structure of this stabilizer subgroup is discussed further in \cref{lemma.stabilizer-subquotient-structure} (following \cite[Propisition VI.2.4]{FarguesScholze.GeometrizationOfTheLocalLanglandsCorrespondence}).
\end{example}

\begin{remark}\label{remark.characteristicp} One can also consider a characteristic $p$ analog: suppose $L=\kappa((\pi))$, where $\kappa$ is a perfect field of characteristic $p$. Then, one can consider the category $\AffPerfd_L$ of affinoid perfectoids over $\Spa(L, \mathcal{O}_L)$, where $\mathcal{O}_L=\kappa[[\pi]]$. The analog of $\mathbb{B}^+_\dR$ in this setting sends $\Spa(R,R^+) \in \AffPerfd_L$ to the completion of $R^+[[t]][1/\pi]$ along $(\xi)$ for $\xi=t-\pi$, which is naturally identified with $R[[\xi]]$; $\mathbb{B}_\dR$ then sends $\Spa(R,R^+)$ to $R((\xi))$. For $G/L$ a connected linear algebraic group one can thus define a $\mathbb{B}^+_\dR$-affine Grassmannian and its Schubert cells as above, and, for $G$ split reductive, \cite[Proposition VI.2.4]{FarguesScholze.GeometrizationOfTheLocalLanglandsCorrespondence} also includes the $v$-analog of \cref{theorem.action-transitive} in this case (\cite[Proposition VI.2.4]{FarguesScholze.GeometrizationOfTheLocalLanglandsCorrespondence} assumes $\kappa$ is a finite field, but this is not necessary for the aspect we consider). Similarly, by working at the edge of the space $\mathcal{Y}$ in the sense of \cite[Remark VI.2.1]{FarguesScholze.GeometrizationOfTheLocalLanglandsCorrespondence}, one obtains a variant with $\mathbb{B}^+_\dR(R)=R[[t]]$, and the analogous construction at the edge of $\mathcal{Y}$ for $L$ a $p$-adic field gives a variant with the Witt vector affine Grassmannian. It is likely possible to obtain a version of \cref{theorem.action-transitive} in these cases. However, as our main motivation is an application in characteristic zero (\cref{ss.application}) and because the key preliminaries for treating non-reductive groups in \cite{HoweKlevdal.AdmissiblePairsAndpAdicHodgeStructuresITranscendenceOfTheDeRhamLattice, HoweKlevdal.AdmissiblePairsAndpAdicHodgeStructuresIITheBiAnalyticAxLindemannTheorem} were developed only in characteristic zero, we do not pursue such generalizations here. 
\end{remark}

\subsection{An application}\label{ss.application}

Our interest in \cref{theorem.action-transitive}  comes mainly from \cite{Howe.InscriptionTwistorsAndPAdicPeriods}, where we defined \emph{inscribed} $\mathbb{B}^+_\dR$-affine Grassmannians and Schubert cells as functors on a larger category of test objects consisting of nilpotent thickenings of $\mathbb{B}^+_\dR$ in order to equip these moduli spaces with natural tangent bundles. In this setting, for computations it is better to define the inscribed Schubert cells directly as $v$-homogeneous spaces, but then one would still like to know that the underlying functor on $\AffPerfd_L$ is the usual affine Grassmannian or Schubert cell. This is implied by \cref{theorem.action-transitive}. Moreover, the theory of \cite{Howe.InscriptionTwistorsAndPAdicPeriods} equips $\Gr_{[\mu]}$ with a Banach-Colmez Tangent Bundle $T_{\Gr_{[\mu]}} / \Gr_{[\mu]}$ (the underlying $v$-sheaf of the inscribed tangent bundle). By constructing the usual tangent bundle of $G$, $T_G/G$, as the restriction of scalars $T_G=\Res_{L[\epsilon]/L} G$ (for $\epsilon^2=0$), \cref{example.ros} below implies (by checking the identity after extending scalars so $[\mu]$ has a representative defined over $L$) that \[ T_{\Gr_{[\mu]}}= \Gr_{[\mu_{T_G}]},  \]
where $[\mu_{T_G}]$ denotes the induced class of cocharacters of $T_G \cong G \ltimes \Lie G$. Note that, for this second computation, even when $G$ is reductive we need the statement of \cref{theorem.action-transitive} for non-reductive groups in order to apply it to the tangent bundle $T_G$, which has non-trivial unipotent radical $\Lie G$.    

\begin{example}\label{example.ros}
Suppose $A$ is a finite $L$-algebra and $H/L$ is a connected linear algebraic group. We write $H_A$ for the base change of $H$ to $A$ and $G=\Res_{A/L} H_A$ for its Weil restriction of scalars back to $L$, which is also a connected linear algebraic group over $L$. Viewing $G$ as a functor on $L$-algebras $B$, we have 
\begin{align*} G(B)&=H_A(A \otimes_L B)=H(A \otimes_L B )=H_B(A \otimes_L B). \end{align*}
In particular, there is a natural map $\iota: H \rightarrow G$ induced by $B \rightarrow A\otimes_L B.$ 

Fix now a cocharacter $\mu$ of $H$ defined over $L$, and let $\mu_G=\iota \circ \mu$. By \cref{remark.homogeneous-space},   $\Gr_{[\mu_G]} / \Spd\, L$ is the \'{e}tale sheafification of 
\[ (R,R^+) \mapsto G(\mathbb{B}^+_\dR(R))/\mathrm{Stab}_{G(\mathbb{B}^+_\dR(R))}(\ast_{\mu_G}). \]
We can rewrite the quotient appearing on the right in terms of $H$ as 
\[ H_{\mathbb{B}^+_\dR(R)}(\mathbb{B}^+_\dR(R) \otimes_L A) / \mathrm{Stab}_{H_{\mathbb{B}^+_\dR(R)}}(\ast_\mu)(\mathbb{B}^+_\dR(R) \otimes_L A) \]
where the stabilizer of $\ast_{\mu}$ can be interpreted here as follows: for $\xi$ a generator of $\mr{Ker}\, \theta \subseteq \mathbb{B}^+_\dR(R)$, and $S$ a $\mathbb{B}^+_\dR(R)$-algebra in which $\xi$ remains invertible, 
\[ \mathrm{Stab}_{H_{\mathbb{B}^+_\dR(R)}}(\ast_\mu)(S) = \xi^{\mu} H(S) \xi^{-\mu} \cap H(S) \leq H(S[1/\xi]). \] 
\end{example}

In the context of the inscribed Schubert cells of \cite{Howe.InscriptionTwistorsAndPAdicPeriods}, it is natural to consider the following question about generalizations of \cref{example.ros}. 

\begin{question}
Let $\Spa(R,R^+)/\Spa\, L$ be affinoid perfectoid and suppose $B$ is a finite locally free $\mathbb{B}^+_\dR(R)$-algebra (the specific case directly relevant to the setup in \cite{Howe.InscriptionTwistorsAndPAdicPeriods} is when furthermore $B^\mr{red}=\mathbb{B}^+_\dR(R)$). For $H/L$ a connected linear algebraic group, we can consider the \'{e}tale sheafification of 
\[\Spa(S,S^+)/\Spa(R,R^+) \mapsto H(\mathbb{B}_\dR(S) \otimes_{\mathbb{B}^+_\dR(R)} B)/ H(\mathbb{B}^+_\dR(S) \otimes_{\mathbb{B}^+_\dR(R)} B).\]
Then, generalizing \cref{example.ros}, for any cocharacter $\mu$ of $H$ we would like to know whether membership in the $H(\mathbb{B}^+_\dR \otimes_{\mathbb{B}^+_\dR(R)} B)$-orbit of $\ast_{\mu}$ in this sheaf can be checked at geometric points. If we write 
\[ G=\Res_{B/\mathbb{B}^+_\dR(R)} H_{B}, \]
then this can be viewed as a question about $G(\mathbb{B}^+_\dR)$ orbits in $\Gr_G/\Spd(R,R^+)$, where the latter is defined to be the \'{e}tale sheafification of
\[ \Spa(S,S^+)/\Spa(R,R^+) \mapsto G(\mathbb{B}_\dR(S))/G(\mathbb{B}^+_\dR(S)). \]
More generally, for any connected linear algebraic group $G/\mathbb{B}^+_\dR(R)$, we can define $\Gr_G/\Spd(R,R^+)$ to be the \'etale quotient $G(\mathbb{B}_\dR)/G(\mathbb{B}^+_\dR)$. Then, for $s_1, s_2 \in \Gr_G(\Spa(R,R^+))$, we ask: if $s_1$ and $s_2$ are in the same $G(\mathbb{B}^+_\dR(C))$-orbit after restricting to any geometric point $\Spa(C,C^+) \rightarrow \Spa(R,R^+)$, is there always an \'{e}tale or $v$-cover $\Spa(S,S^+) \rightarrow \Spa(R,R^+)$ such that $s_1$ and $s_2$ are in the same $G(\mathbb{B}^+_\dR(S))$-orbit? If not, what further conditions are necessary? 
\end{question}

\subsection{Outline} In \cref{s.preliminaries} we recall and develop some preliminaries on the $\mathbb{B}^+_\dR$-affine Grassmannian and its Schubert cells. In \cref{s.proof} we prove \cref{theorem.action-transitive}.

\subsection{Acknowledgements} During the preparation of this work, the author was supported by the NSF through grants DMS-2201112 and DMS-2501816. We thank Peter Scholze for providing the reference  \cite[Proposition VI.2.4]{FarguesScholze.GeometrizationOfTheLocalLanglandsCorrespondence} and for helpful correspondence related to \cref{prop.closed-schubert-functor}. We thank Kestutis \v{C}esnavi\v{c}ius, Christian Klevdal, and Alex Youcis for helpful conversations. 

\section{Preliminaries}\label{s.preliminaries}

\subsection{Perfectoids and period sheaves}
Let $L$ be a $p$-adic field, that is, a complete discretely valued extension of $\mathbb{Q}_p$ with perfect residue field. We consider the category $\AffPerfd_L$ of affinoid perfectoid spaces over $L$. There are natural analytic, \'{e}tale, and $v$-topologies on $\AffPerfd_L$ defined as in \cite[Top of p.3]{Scholze.EtaleCohomologyOfDiamonds}.

We recall (from, e.g., \cite[Lecture 15]{ScholzeWeinstein.BerkeleyLecturesOnPAdicGeometryAMS207}, \cite[\S1]{CesnaviciusYoucis.TheAnalyticTopologySufficesForTheBplusDrGrassmannian} or \cite[Section 2.3]{HoweKlevdal.AdmissiblePairsAndpAdicHodgeStructuresIITheBiAnalyticAxLindemannTheorem}) that we have Fontaine functors $R \mapsto \mathbb{B}^+_\dR(R)$ and $R \mapsto \mathbb{B}_\dR(R)$ from perfectoid $L$-algebras to $L$-algebras, and a functorial surjection $\theta: \mathbb{B}^+_\dR(R) \rightarrow R$, whose kernel is a principal ideal. For $\xi$ any generator of $\mr{Ker}\, \theta$, $\xi$ is not a zero divisor and $\mathbb{B}_\dR(R)=\mathbb{B}^+_\dR[1/\xi]$. We equip $\mathbb{B}^+_\dR(R)$ with the $\mr{Ker}\, \theta$-adic filtration and $\mathbb{B}_\dR(R)$ with the induced filtration. Both are complete and separated; in particular, for $K/L$ a perfectoid field, $\mathbb{B}^+_\dR(K)$ is a complete discrete valuation ring over $L$ with residue field $K$ and fraction field $\mathbb{B}_\dR(K)$. Moreover, the assignments 
\[ \Spa(R,R^+) \mapsto \mathbb{B}^+_\dR(R) \textrm{ and } \Spa(R,R^+) \mapsto \mathbb{B}_\dR(R) \]
are sheaves for the $v$-topology on $\AffPerfd_L$. We will also consider the structure $v$-sheaf $\mathcal{O}$ on $\AffPerfd_L$ which sends $\Spa(R,R^+)$ to $R$; note that $\mathcal{O}=\mathbb{B}^+_\dR /\Fil^1\mathbb{B}^+_\dR$ (where the quotient can be taken even as presheaves). 

\subsection{The $\mathbb{B}^+_\dR$-affine Grassmannian}

\begin{definition} Let $L$ be a $p$-adic field. For $G/L$ a connected linear algebraic group, we consider the presheaf on $\AffPerfd_L$
\[ \Gr_G^{\mr{pre}}: \Spa(R,R^+) \mapsto G(\mathbb{B}_\dR(R))/G(\mathbb{B}^+_\dR(R)). \]
The $\mathbb{B}^+_\dR$-affine Grassmannian $\Gr_G$ is the \'{e}tale sheafification of $\Gr_G^{\mr{pre}}$. 
\end{definition}

We recall some properties of $\Gr_G$ (essentially from \cite{FarguesScholze.GeometrizationOfTheLocalLanglandsCorrespondence, ScholzeWeinstein.BerkeleyLecturesOnPAdicGeometryAMS207, CesnaviciusYoucis.TheAnalyticTopologySufficesForTheBplusDrGrassmannian}; see also \cite[\S5]{HoweKlevdal.AdmissiblePairsAndpAdicHodgeStructuresIITheBiAnalyticAxLindemannTheorem}).

\begin{proposition}\label{prop.aff-grass-properties} Let $L$ be a $p$-adic field and let $G/L$ be a connected linear algebraic group. 
\begin{enumerate}
    \item Consider the presheaf $\mathcal{F}$ on $\AffPerfd_L$ sending $\Spa(R,R^+)$ to
    \[ \left\{ (\mc{E}, \iota)\, | \textrm{ $\mc{E}$ a $G$-torsor on $(\Spec \mathbb{B}^+_\dR(R))_\et$ and $\iota \in \mc{E}(\mathbb{B}_\dR(R))$} \right\}/\sim. \] 
    Then, writing $G_{\mathbb{B}^+_\dR(R)}$ for the trivial $G$-torsor, the map
    \[ G(\mathbb{B}_\dR) \rightarrow \mc{F},\; g \in G(\mathbb{B}_\dR(R)) \mapsto (G_{\mathbb{B}^+_\dR(R)}, g) \]
    factors uniquely through an isomorphism $\Gr_G \xrightarrow{\sim} \mc{F}$. 
    \item $\Gr_G$ is a sheaf for the $v$-topology on $\AffPerfd_L$.
    \item $\Gr_G$ is the sheafification of $\Gr_G^{\mr{pre}}$ for the analytic topology. 
\end{enumerate} 
\end{proposition}
\begin{proof}
For $G$ reductive, parts (1) and (2) were treated in \cite[Proposition 19.1.2]{ScholzeWeinstein.BerkeleyLecturesOnPAdicGeometryAMS207}, except that in \cite[Definition 19.1.1]{ScholzeWeinstein.BerkeleyLecturesOnPAdicGeometryAMS207} $G$ is taken over a complete algebraically closed extension $C/\mathbb{Q}_p$ and implicitly base-changed along a choice of a section $C \rightarrow \mathbb{B}^+_\dR(C)$ of $\theta:\mathbb{B}^+_\dR(C) \rightarrow C$ so that, in particular, one is working over $\Spd\, C$. Parts (1) and (2) were then treated in \cite[Proposition VI.1.9]{FarguesScholze.GeometrizationOfTheLocalLanglandsCorrespondence} for a more general notion of $\mathbb{B}^+_\dR$ associated to a general divisor on the Fargues-Fontaine curve (or the related space $\mathcal{Y}$), but still under the assumption of $G$ reductive and now with $L/\mathbb{Q}_p$ a finite extension. These proofs used the Tannakian formalism and a reduction to the $v$-stack property for the assignment sending $\Spa(R,R^+)$ to the category of finite projective $R$-modules as established in \cite[Lemma 17.1.8]{ScholzeWeinstein.BerkeleyLecturesOnPAdicGeometryAMS207}\footnote{Note that the proof of \cite[Lemma 17.1.8]{ScholzeWeinstein.BerkeleyLecturesOnPAdicGeometryAMS207} is written for finite projective $R$-modules but its statement combines this $v$-stack property with the equivalence between finite projective $R$-modules and locally free $\mathcal{O}_{\Spa(R,R^+)}$-modules of finite rank of \cite[Theorem 2.7.7]{KedlayaLiu.RelativepAdicHodgeTheoryFoundations}.}; in particular, these proofs do not use the reductive hypothesis and apply also in our setting where we work over an arbitrary $p$-adic field $L$ (on $\AffPerfd_L$, since $L$ is a $p$-adic field $\mathbb{B}^+_\dR$ is canonically valued in $L$-algebras thus there is no implicit choice of a section needed to make the base change as in \cite[Definition 19.1.1]{ScholzeWeinstein.BerkeleyLecturesOnPAdicGeometryAMS207}). A more detailed and explicit proof of the statements of (1) and (2) using many of the same tools appears in \cite[proof of Proposition 2.1]{CesnaviciusYoucis.TheAnalyticTopologySufficesForTheBplusDrGrassmannian} --- note that the statement of \cite[Proposition 2.1]{CesnaviciusYoucis.TheAnalyticTopologySufficesForTheBplusDrGrassmannian} is given without the reductive hypothesis but with the hypothesis that $L/\mathbb{Q}_p$ is a finite extension, however, this latter hypothesis is not needed for the proof. 

For part (3), \cite[Theorem 3.1]{CesnaviciusYoucis.TheAnalyticTopologySufficesForTheBplusDrGrassmannian} gives the reductive case (with the same caveat that it is stated for $L/\mathbb{Q}_p$ a finite extension but the proof works for any $p$-adic field $L$), but, since we are working in characteristic zero, the result extends to the non-reductive case using the following argument that was suggested to us by Kestutis \v{C}esnavi\v{c}ius: We need to show that, for $\Spa(R,R^+) \in \AffPerfd_L$  and $\mc{E}$ a $G$-torsor on $(\Spec \mathbb{B}^+_\dR(R))_\et$, there is a cover of $\Spa(R,R^+)$ by rational opens $\Spa(S,S^+)$ such that the pullback of $\mc{E}$ along $\Spec \mathbb{B}^+_\dR(S) \rightarrow \Spec \mathbb{B}^+_\dR(R)$ is a trivial $G$-torsor on $(\Spec \mathbb{B}^+_\dR(S))_\et$. We write $U$ for the unipotent radical of $G$ and $M:=G/U$. Since we know the reductive case, by replacing $\Spa(R,R^+)$ with a covering by rational opens, we may assume the push-out of $\mc{E}$ to an $M$-torsor, $\overline{\mc{E}}=\mc{E} \times^G M$, is trivial. Fixing such a trivialization $\overline{s} \in \overline{\mc{E}}(\mbb{B}^+_\dR(R))$, the lifts of $\overline{s}$ to $\mc{E}$ form a $U$-torsor over $\Spec \mbb{B}^+_\dR(R)$. Since $U$ is a unipotent group and $L$ is of characteristic zero, $U$ is an iterated extension of vector groups, and thus we deduce $H^1( (\Spec \mathbb{B}^+_\dR(R))_\et, U)=0$ (from the vanishing of \'{e}tale cohomology for vector groups on an affine scheme). Thus, any $U$-torsor on $\Spec \mathbb{B}^+_\dR(R)$ is trivial, so there is a lift $s$ of the trivialization $\overline{s}$ already over $\Spec \mbb{B}^+_\dR(R)$. Alternatively, by \cite[Remark 3.2]{CesnaviciusYoucis.TheAnalyticTopologySufficesForTheBplusDrGrassmannian}, it suffices to make this argument just for $R/L$ a perfectoid field.
\end{proof}

\subsection{Schubert cells}\label{ss.schubert-cells}
Let $L$ be a $p$-adic field, let $\overline{L}$ be an algebraic closure of $L$, let $G/L$ be a connected linear algebraic group, and let $[\mu]$ be a conjugacy class of cocharacters of $G_{\overline{L}}$. We write $L([\mu]) \subseteq \overline{L}$ for the reflex field, i.e. the field of definition of $[\mu]$, i.e. the fixed field for the stabilizer of $[\mu]$ for the action of $\Gal(\overline{L}/L)$ on conjugacy classes of cocharacters of $G_{\overline{L}}$. The field $L([\mu])$ is a finite extension of $L$ and thus also a $p$-adic field. 

Viewing $\Gr_G \times_{\Spd L} \Spd L([\mu])$ as a $v$-sheaf on $\AffPerfd_{L([\mu])}$, we define the Schubert cell $\Gr_{[\mu]}$ to be the subfunctor whose values on   $\Spa(R,R^+)$ are those sections whose pullback along any geometric point $\Spa(C,C^+) \rightarrow \Spa(R,R^+)$ lies in
\begin{multline*} G(\mathbb{B}^+_\dR(C)) \xi^{\mu} G(\mathbb{B}^+_\dR(C)) / G(\mathbb{B}^+_\dR(C)) \subseteq \\G(\mathbb{B}_\dR)(C)/G(\mathbb{B}^+_\dR(C))= \Gr_G(\Spa(C,C^+))\end{multline*}
for one (equivalently, any) choice of a representative $\mu \in [\mu]$, an embedding $\overline{L} \rightarrow C$ extending the given $L([\mu])\rightarrow C$, and generator $\xi$ of $\mathrm{Ker}\,\theta \subseteq \mathbb{B}^+_\dR(C)$. It is clear from the $v$-sheaf property of $\Gr_G$ and this definition that $\Gr_{[\mu]}$ is also a $v$-sheaf.

\subsubsection{}\label{sss.closed-faithful}
To prove \cref{theorem.action-transitive} in the reductive case following \cite[Proof of Proposition VI.2.4]{FarguesScholze.GeometrizationOfTheLocalLanglandsCorrespondence}, we need to know that faithful representations induce closed immersion of Schubert cells in order to carry out a reduction to $\GL_n$. This is \cref{prop.closed-schubert-functor} below --- note that in \cite[Proof of Proposition VI.2.4]{FarguesScholze.GeometrizationOfTheLocalLanglandsCorrespondence} it is claimed as a consequence of \cite[Proposition 20.3.7/19.2.3]{ScholzeWeinstein.BerkeleyLecturesOnPAdicGeometryAMS207}, which gives a closed immersion on the full affine Grassmannian, but a small supplementary argument is required (see \cref{remark.necessity-of-arg}). 

\subsubsection{}In the proof of \cref{prop.closed-schubert-functor}, we will need the Bruhat order on geometric conjugacy classes of cocharacters: we say $[\nu] \leq [\mu]$ if, after fixing a maximal torus $T\subseteq G$ and Borel subgroup $T \subseteq B \subseteq G$, for $\mu$ and $\nu$ the associated dominant representatives, $\mu - \nu$ is a sum of simple coroots. We write $\Gr_{\leq[\mu]} \subseteq \Gr_G \times_{\Spd L} \Spd L([\mu])$ for the Schubert variety defined to consist of those sections that, at all geometric points, lie in the union of the Schubert cells $\Gr_{[\nu]}$ for $[\nu] \leq [\mu]$ in the Bruhat order. For $G$ reductive, $\Gr_{\leq [\mu]}$ is closed in $\Gr_G$ by \cite[Proposition 19.2.3]{ScholzeWeinstein.BerkeleyLecturesOnPAdicGeometryAMS207}. 

\begin{proposition}\label{prop.closed-schubert-functor}
    Let $G/L$ be a connected reductive group and let $\rho: G \rightarrow \GL_n$ be a faithful representation. Then, for any conjugacy class $[\mu]$ of cocharacters of $G_{\overline{L}}$, the map $\Gr_{[\mu]} \rightarrow \Gr_{[\rho \circ \mu]} \times_{\Spd L} \Spd L([\mu])$
    induced by $\rho$ is a closed immersion. 
\end{proposition}
\begin{proof}
    We may extend scalars to assume $L=L([\mu])$. By \cite[Lemma 19.1.5]{ScholzeWeinstein.BerkeleyLecturesOnPAdicGeometryAMS207}, the map $\Gr_G \rightarrow \Gr_{\GL_n}$ induced by $\rho$ is a closed immersion. As noted above, $\Gr_{\leq [\mu]}$ is closed in $\Gr_G$, so it is also closed in $\Gr_{\GL_n}$. Thus, it suffices to show that $\Gr_{[\mu]} = \Gr_{[\rho \circ \mu]} \times_{\Gr_{\GL_n}} \Gr_{\leq [\mu]}.$
    Since $\Gr_{[\mu]}$ is contained in this intersection, we need only show that $[\nu] < [\mu]$ implies $[\rho \circ \nu] \neq [\rho \circ \mu]$. This follows from \cref{lemma.cocharacter-pushforward}.
\end{proof}

\begin{remark}\label{remark.necessity-of-arg}
    We note that typically $\Gr_{[\mu]}\neq\Gr_{[\rho \circ \mu]} \times_{\Gr_{\GL_n}} \Gr_{G}$. For example, consider the case of $G$ a maximal torus in $\GL_n$ --- in this case $[\mu]=\{\mu\}$ has a unique element and the intersection also contains $\Gr_{\{w \cdot \mu\}}$ for each element $w$ of the Weyl group of $\GL_n$, the symmetric group $S_n$. Thus some further argument as above is necessary to deduce \cref{prop.closed-schubert-functor} from \cite[Lemma 19.1.5]{ScholzeWeinstein.BerkeleyLecturesOnPAdicGeometryAMS207}.
\end{remark}

\subsubsection{}\label{sss.group-theoretic-lemma-and-conj-class-intro}
We used the following group theoretic lemma in the proof of \cref{prop.closed-schubert-functor}. To state it, we recall that, for $V$ an $n$-dimensional vector space over $L$, the  conjugacy classes $[\mu]$ of cocharacters of $\GL(V)$ are parameterized by tuples of integers $k_1 \geq k_2 \geq \ldots \geq k_n$ --- the conjugacy class $[\mu]$ corresponds to the tuple such that there exists a basis of $V$ where $\mathbb{G}_m$ acts through $\mu$ on the $i$th basis vector by the character $t\mapsto t^{k_i}$. In this case, in addition to the Bruhat order, we also have the dictionary order, where $(k_1, \ldots, k_n) > (k_1', \ldots, k_n')$ if the first entry $i$ where they differ has $k_i > k_i'$. In this indexing, $(k_1, \ldots, k_n) > (k_1', \ldots, k_n')$ in the Bruhat order if $(k_1, \ldots, k_n)$ can be obtained from $(k_1', \ldots, k_n')$ by adding vectors of the form $s_i^\vee=(0,\ldots, 0, \overbrace{1,-1}^{i,i+1}, 0, \ldots,0)$, $1\leq i < n$, so the dictionary order refines the Bruhat order (i.e. if $[\mu'] < [\mu]$ in Bruhat order then also $[\mu'] < [\mu]$ in dictionary order).  

\begin{lemma}\label{lemma.cocharacter-pushforward}
    Let $G/L$ be a connected reductive group and let $[\mu]$ be a conjugacy class of cocharacters of $G$. Then, for any faithful representation $\rho: G \rightarrow \GL_n$, if $[\nu]<[\mu]$ in the Bruhat order then $[\rho \circ \nu] < [\rho \circ \mu]$ in the dictionary order.  
\end{lemma}
\begin{proof}
We extend scalars everywhere to $\overline{L}$, then fix a maximal torus and Borel $T \leq B \leq G$. Let $\nu$ and $\mu$ be the associated dominant representatives. By applying 
\cite[Lemma 2.3]{Rapoport.APositivityPropertyOfTheSatakeIsomorphism} repeatedly, we find that $[\nu] < [\mu]$ implies that we can find simple coroots $\alpha_1^\vee, \ldots, \alpha_k^\vee$ such that 
$\mu=\nu+\alpha_1^\vee + \ldots + \alpha_k^\vee$ and furthermore $\nu+\alpha_1^\vee + \ldots + \alpha_i^\vee$ is also dominant for any $1 < i < k$. In particular, it suffices to treat the case where $\mu=\nu+\alpha^\vee$ for a simple coroot $\alpha^\vee$. 

Writting $L=\langle \ell \rangle$ for the line associated to $\alpha$ in $\Lie G$, we can view $L$ as a line inside of $\Lie \GL_n$. The generator $\ell$ is nilpotent, and by conjugating the embedding into $\GL_n$ we can assume it is in Jordan form, so a product of standard nilpotent Jordan blocks. The cocharacter $\nu$ factors through the stabilizer $H \leq \GL_n$ of $L$, and lies in a maximal torus $T'$ of $H$; by conjugating by an element of $H$ we can assume this maximal torus is the intersection of the diagonal matrices in $\GL_n$ with $H$--- this contains the group of diagonal matrices $T_0$ that are constant on each block and $T'/T_0$ is equal to the diagonal torus in the principal $\mathrm{PSL}_2$ generated by $L$ (here we mean either the principal $\mathrm{PSL}_2$ or the quotient of the principal $\SL_2$ by its center, depending on whether there is a Jordan block of even dimension or not). In particular, there is a non-negative even integer $k$ such that the weights of $\rho \circ \nu$ on each Jordan block are of the form $m, m-k, m-2k,\ldots$ for some integer $m$ (where $m$ depends on the block). Adding $\alpha^\vee$ to $\nu$ to obtain $\mu$ changes the weights on the blocks by replacing $k$ with $k+2$ and each $m$ with $m + r -1$, where $r$ is the size of the given block. In particular, it follows that $[\rho \circ \mu]=[\rho\circ(\nu+\alpha^\vee)]$, viewed as a cocharacter of $\GL_n$, is larger than $[\rho \circ \nu]$ in the dictionary order. Indeed, if $\rho \circ \nu$ has weights $k_1 \geq \ldots \geq k_n$, then, for $M$ the maximum of $m+r-1$ over the Jordan blocks of size $r>1$, the weights of $\rho \circ \mu$ will be $k_1 \geq \ldots \geq k_i \geq M \geq \ldots $ where $i$ is the last index where $k_i \geq M$ that can be taken from a Jordan block of size $1$ (we set $i=0$ if no such index exists), and one checks immediately that $M > k_{i+1}$. 
\end{proof}

\subsection{Canonical points and their stabilizers}\label{ss.canonical-point-and-stabilizer}
Let $L$ be a $p$-adic field and let $G/L$ be a connected linear algebraic group. For $\mu: \mathbb{G}_m \rightarrow G$ a cocharacter, there is a canonical point $\ast_{\mu}: \Spd L \rightarrow \Gr_{[\mu]}$
that, on $\Spa(R,R^+)$, is given by $\xi^\mu G(\mathbb{B}^+_\dR(R))$ for any generator $\xi$ of $\mathrm{Ker}\,\theta$.

We now describe the stabilizer of $\ast_{\mu}$ in $G(\mathbb{B}^+_\dR)$ via its subquotients for the filtration by intersections with the principal congruence subgroups
\begin{equation}\label{eq.def-pcsg} G_i:= \mathrm{Ker} \left( G(\mathbb{B}^+_\dR) \rightarrow G(\mathbb{B}^+_\dR/\Fil^i \mathbb{B}^+_\dR) \right).\end{equation}
Before stating the result, we note that for $\nu: \mathbb{G}_m \rightarrow G$ a cocharacter and $W$ a representation of $G$, there is a descending filtration on $W$ given by setting $\Fil^i_{\nu}(W)$ to be the sum of the weight subspaces where $\mathbb{G}_m$ acts through $\nu$ by $t \mapsto t^j$, $j \geq i$. We write $P_{\nu} \subseteq G$ for the parabolic subgroup preserving this filtration (for all $W$). In the following, we view $\Lie G$ as a representation via the adjoint action of $G$ and we write $\{i\}$ for a Breuil-Kisin-Fargues twist, i.e. for tensor product with  $\mathcal{O}\{i\}:=\gr^{i}\mathbb{B}_\dR$.  

\begin{lemma}\label{lemma.stabilizer-subquotient-structure}Let $L$ be a $p$-adic field and let $G/L$ be a connected linear algebraic group. For $\mu: \mathbb{G}_m \rightarrow G$ a cocharacter, let $H=\mathrm{Stab}_{G(\mathbb{B}^+_\dR)}(\ast_{\mu})$ and, for $i \geq 0$, let $H_i=H \cap G_i$. There are canonical and functorial identifications of $v$-sheaves
    \[ H_0/H_1 = P_{\mu^{-1}}(\mathcal{O}) \subseteq G(\mathcal{O})=G_0/G_1 \]
and, for $i\geq 1$, 
    \[ H_i/H_{i+1}=\Fil^{-i}_{\mu^{-1}}(\Lie G)\{i\} \subseteq (\Lie G)\{i\}=G_i/G_{i+1}.\]
\end{lemma}
\begin{proof}
We follow the second and third paragraphs of \cite[Proof of Proposition VI.2.4]{FarguesScholze.GeometrizationOfTheLocalLanglandsCorrespondence}. Thus, we fix a closed embedding $\rho: G \hookrightarrow \GL_n$; conjugating, we may assume $\rho \circ \mu$ sends $t$ to the diagonal matrix $\mr{diag}(t^{k_1},\ldots, t^{k_n})$, $k_1 \geq k_2 \geq \ldots \geq k_n$. The claims for $G=\GL_n$ with cocharacter $\rho \circ \mu$ are then readily verified by writing, for any $\Spa(R,R^+) \in \AffPerfd_L$ and generator $\xi$ of $\mr{Ker}\,\theta$,  $\mr{Stab}_{\GL_n(\mathbb{B}^+_\dR(R))}(\ast_{\rho \circ \mu})$ as 
\[ \mr{diag}(\xi^{k_1},\ldots, \xi^{k_n}) \GL_n(\mathbb{B}^+_\dR(R))  \mr{diag}(\xi^{-k_1},\ldots, \xi^{-k_n}) \cap \GL_n(\mathbb{B}^+_\dR(R)), \]
which is the group of matrices $A=(A_{ij}) \in \GL_n(\mathbb{B}^+_\dR(R))$ with $A_{ij} \in \xi^{k_i - k_j} \mathbb{B}^+_\dR(R)$. 

For our general $G$, it then follows that 
\begin{equation}\label{eq.inclusion-i-0} H_0/H_1 \subseteq P_{\mu^{-1}}(\mathcal{O})=P_{\rho \circ \mu^{-1}}(\mathcal{O}) \cap G(\mathcal{O})\end{equation}
and that, for $i \geq 1$, 
\begin{equation}\label{eq.inclusion-i-g-0} H_i / H_{i+1} \subseteq \Fil_{\mu^{-1}}^{-i}(\Lie G)\{i\}= \Lie G\{i\} \cap\Fil_{\rho \circ \mu^{-1}}^{-i}(\Lie \GL_n)\{i\}.\end{equation}
On the other hand, we have $P_{\mu^{-1}}(\mathbb{B}^+_\dR) \subseteq H_0$ and, for each root $\alpha$ appearing in the unipotent radical of the opposite parabolic $P_{\mu}$, if we write the associated one-parameter subgroup of $G$ as $U_{\alpha}$ and fix an isomorphism $U_{\alpha} \cong \mathbb{G}_a$, then $\Fil^{\langle \alpha, \mu \rangle}\mathbb{B}^+_\dR \subseteq \mathbb{B}^+_\dR = U_{\alpha}(\mathbb{B}^+_\dR)$ is also contained in $H_0$. Together the existence of these subgroups shows that both of the inclusions \cref{eq.inclusion-i-0} and \cref{eq.inclusion-i-g-0} are in fact equalities. 
\end{proof}

As a first consequence of \cref{lemma.stabilizer-subquotient-structure}, we obtain the following result, which we will use to promote the $v$-topology version of \cref{theorem.action-transitive} in the reductive case to the desired statement in the \'{e}tale topology. 

\begin{proposition}\label{prop.v-quotient-is-etale-quotient}
    Let $L$ be a $p$-adic field and let $G/L$ be a connected linear algebraic group. For $\mu: \mathbb{G}_m \rightarrow G$ a cocharacter, let $H=\mathrm{Stab}_{G(\mathbb{B}^+_\dR)}(\ast_{\mu})$. Then, the \'{e}tale quotient $G(\mathbb{B}^+_\dR)/H$ agrees with the $v$-quotient $G(\mathbb{B}^+_\dR)/H$  on $\AffPerfd_L$. 
\end{proposition}
\begin{proof}
 It suffices to show that the map $G(\mathbb{B}^+_\dR) \rightarrow G(\mathbb{B}^+_\dR)/H$, where the quotient on the right is formed in the $v$-topology, 
is surjective already in the \'{e}tale topology. 

We write $G_i$ for the principal congruence subgroups as in \cref{eq.def-pcsg} and $H_i=G_i \cap H$ as in \cref{lemma.stabilizer-subquotient-structure}.  
To establish the surjectivity, we first note that, for any $n\geq 0$, the map $G(\mathbb{B}^+_\dR)=G_0 \rightarrow G_0/G_n$ is surjective already as a map of pre-sheaves. Indeed, $G$ is smooth affine and, for any $(R,R^+) \in \AffPerfd_L$,
\[ \mathbb{B}^+_\dR(R) \twoheadrightarrow \mathbb{B}^+_\dR(R)/\Fil^n \mathbb{B}^+_\dR(R) \]
is the limit of the nilpotent thickenings \[ \mathbb{B}^+_\dR(R)/\Fil^{n+i} \mathbb{B}^+_\dR(R) \twoheadrightarrow \mathbb{B}^+_\dR(R)/\Fil^n \mathbb{B}^+_\dR(R),\]
so this follows from the lifting property for smooth algebras \cite[\href{https://stacks.math.columbia.edu/tag/07K4}{Tag 07K4}]{stacks-project}. 

Now, if we fix $n$ large enough that $\Fil^{-n}_{\mu^{-1}}(\Lie G)=\Lie G$, then \cref{lemma.stabilizer-subquotient-structure} implies that $G_n \leq H_0=H$ and thus $G_0 \rightarrow G_0/H_0$ factors through $G_0/G_n$. Since $G_0 \rightarrow G_0/G_n$ is surjective already as a map of presheaves, it suffices to show that the $v$-surjection
$G_0/G_n \rightarrow G_0/H_0$
is in fact surjective in the \'{e}tale topology. To that end, note that we can factor this map into a sequence of surjections of $v$-sheaves
\[ G_0/H_0 \twoheadleftarrow G_0/H_1 \twoheadleftarrow \ldots \twoheadleftarrow G_0/H_n = G_0/G_n. \]
For each $0 \leq i < n$, $G_0/H_{i+1} \rightarrow G_0/H_{i}$ is a torsor for the $v$-sheaf of groups $H_i/H_{i+1}$. By \cref{lemma.stabilizer-subquotient-structure},  $H_0/H_{1}=P_{\mu^{-1}}(\mc{O})$, thus, after pullback to $\Spa(R,R^+)\in \AffPerfd_L$, \cite[Theorem 19.5.2]{ScholzeWeinstein.BerkeleyLecturesOnPAdicGeometryAMS207} implies the $i=0$ surjection is split after passing to an \'{e}tale cover $\Spa(S,S^+) \rightarrow \Spa(R,R^+)$. When $i>0$, \cref{lemma.stabilizer-subquotient-structure} implies that the pullback of $H_i/H_{i+1}$ to $\Spa(S,S^+)$ is a vector group isomorphic to $\mathcal{O}^m$ for some $m \geq 0$, thus the $i>0$ surjections are split without passing to any further cover using the vanishing of $H^1(\Spa(S,S^+)_v, \mc{O})$ established in \cite[Theorem 17.1.3]{ScholzeWeinstein.BerkeleyLecturesOnPAdicGeometryAMS207}. 

\end{proof}

\section{Proof of the main result}\label{s.proof}
In this section we prove \cref{theorem.action-transitive}. We split up the proof into three steps: first, in \cref{ss.glV} we treat the case of $G=\GL(V)$. Following identically the argument of \cite[Proof of Proposition VI.2.4]{FarguesScholze.GeometrizationOfTheLocalLanglandsCorrespondence}, in this case we obtain transitivity even in the analytic topology. In \cref{ss.red-case} we then follow the argument of \cite[Proof of Proposition VI.2.4]{FarguesScholze.GeometrizationOfTheLocalLanglandsCorrespondence} to obtain transitivity for the general reductive case in the $v$-topology (invoking \cref{prop.closed-schubert-functor} to fill in a small detail of the proof as explained in \cref{sss.closed-faithful}). Combined with \cref{prop.v-quotient-is-etale-quotient}, this gives \cref{theorem.action-transitive} in the reductive case. In \cref{ss.general-case} we finish by explaininig why the proof of the reductive case does not immediately generalize to all connected linear algebraic groups and then giving an argument similar to those used in \cite[\S5.2]{HoweKlevdal.AdmissiblePairsAndpAdicHodgeStructuresIITheBiAnalyticAxLindemannTheorem} that instead deduces the general case of \cref{theorem.action-transitive} from the reductive case.

\subsection{The case of $\GL(V)$}\label{ss.glV}

\begin{proposition}\label{prop.gln}
Let $L$ be a $p$-adic field and let $G=\GL(V)$ for $V$ a finite dimensional $L$-vector space. For any $[\mu]$,     the action of $G(\mathbb{B}^+_\dR)$ on the Schubert cell $\Gr_{[\mu]}$ is transitive in the analytic topology on $\AffPerfd_{L}$.
\end{proposition}
\begin{proof}
We follow the argument given in \cite[First paragraph of the proof of proposition VI.2.4]{FarguesScholze.GeometrizationOfTheLocalLanglandsCorrespondence}: 
We write $k_1 \geq \ldots \geq k_n$ for the non-decreasing list of integers associated to $[\mu]$ as in \cref{sss.group-theoretic-lemma-and-conj-class-intro}. 
    Using the moduli interpretation of \cref{prop.aff-grass-properties}-(1) and pushing out a $\GL(V)$-torsor along the standard representation on $V$, giving a point $s \in \Gr_G(\Spa(R,R^+))$ is equivalent to giving a finite projective $\mathbb{B}^+_\dR(R)$-module $M$ and a trivialization $M \otimes_{\mathbb{B}^+_\dR(R)} \mathbb{B}_\dR(R) = V \otimes_L \mathbb{B}_\dR(R)$; i.e. to giving a  $\mathbb{B}^+_\dR(R)$-lattice $M$ in $V \otimes_L \mathbb{B}_\dR(R)$. Fixing a generator $\xi$ for $\mr{Ker}\,\theta \subseteq \mathbb{B}^+_\dR(R)$, in this intepretation, $s \in \Gr_{[\mu]}(\Spa(R,R^+))$ if, for any geometric point $\Spa(C,C^+) \rightarrow \Spa(R,R^+)$, we can choose a basis $f_1, \ldots, f_n$ of $V \otimes_L \mathbb{B}^+_\dR(C)$ so that 
    \[ M \otimes_{\mathbb{B}^+_\dR(R)} \mathbb{B}^+_\dR(C) \subseteq V \otimes_{L} \mathbb{B}_\dR(C) \]
    is spanned by $\xi^{k_1} f_1, \xi^{k_2} f_2, \ldots, \xi^{k_n} f_n$. We need to show that there is then also an analytic cover $\Spa(S,S^+) \rightarrow \Spa(R,R^+)$ and a basis $f_1, \ldots, f_n$ of $V \otimes_L \mathbb{B}^+_\dR(S)$ such that $\xi^{k_1} f_1, \xi^{k_2} f_2, \ldots, \xi^{k_n} f_n$ is a basis for $M \otimes_{\mathbb{B}^+_\dR(R)} \mathbb{B}^+_\dR(S)$. 
    
    To that end, by \cite[Proposition 19.4.2 or its proof]{ScholzeWeinstein.BerkeleyLecturesOnPAdicGeometryAMS207}, note that for $i \in \mathbb{Z}$,
    \[ \Fil^i (V \otimes_L R) = \left(\xi^i M \cap V \otimes_L \mathbb{B}^+_\dR(R) \right)/\left( \xi^{i} M \cap \xi (V \otimes_L \mathbb{B}^+_\dR(R))\right) \]
    is a local direct summand of $V \otimes_L R$ of rank $m_i:=\#\{j,\,|\, 1\leq j \leq n \textrm{ and } k_j \leq -i\}$ 
    and its formation is stable under the change of base corresponding to any map $\Spa(S,S^+)\rightarrow \Spa(R,R^+)$ in $\AffPerfd_{L}$.  
    Thus, by \cite[Theorem 2.7.7]{KedlayaLiu.RelativepAdicHodgeTheoryFoundations}, we may replace $\Spa(R,R^+)$ with a disjoint union of rational opens $\Spa(S,S^+)$ to assume that each $\Fil^i (V \otimes_L R)/\Fil^{i+1}(V \otimes_L R)$ is free. We may then choose a basis $e_1, \ldots, e_n$ of $V \otimes_L R$ such that, for each $i$, $e_1, \ldots, e_{m_i}$ is a basis for $\Fil^i (V \otimes_L R)$. Then, for each $1 \leq j \leq n$, we choose a lift of $e_{j}$ to $f_{j} \in \xi^{-k_j} M \cap V \otimes_L \mathbb{B}^+_\dR(R)$, and this is the desired basis for $V \otimes_L \mathbb{B}^+_\dR(R)$.
\end{proof}

\subsection{The case of $G$ reductive}\label{ss.red-case}

\begin{proposition}\label{prop.red-case}
Let $L$ be a $p$-adic field, let $\overline{L}$ be an algebraic closure, let $G/L$ be a connected reductive group, and let $[\mu]$ be a conjugacy class of cocharacters of $G_{\overline{L}}$ with reflex field $L([\mu])$. 
    The action of $G(\mathbb{B}^+_\dR)$ on the Schubert cell $\Gr_{[\mu]}$ is transitive in the \'{e}tale topology on $\AffPerfd_{L([\mu])}$.
\end{proposition}
\begin{proof}
By base change to a finite extension $L'/L$, it suffices to assume we have a representative $\mu$ defined over $L$. We write $\ast_\mu$ for the associated point as in \cref{ss.canonical-point-and-stabilizer}. Invoking \cref{prop.v-quotient-is-etale-quotient}, it suffices to show that the orbit map for the action of $G(\mathbb{B}^+_\dR)$ on $\ast_\mu$ induces an isomorphism of $v$-sheaves
\[ G(\mathbb{B}^+_\dR)/\mathrm{Stab}_{G(\mathbb{B}^+_\dR)}(\ast_\mu) = \Gr_{[\mu]}. \]
To show this map is an isomorphism of $v$-sheaves, we follow the argument given in \cite[Fourth paragraph of the proof of Proposition VI.2.4]{FarguesScholze.GeometrizationOfTheLocalLanglandsCorrespondence}. 

    To that end, we fix a faithful representation $\rho:G \rightarrow \GL_n$ and let $\mu'=\rho\circ \mu$ denote the induced character of $\GL_n$. By \cref{prop.gln}, the orbit map for $\ast_{\mu'}$ induces an isomorphism \[ \GL_n(\mathbb{B}^+_\dR)/\mathrm{Stab}_{\GL_n(\mathbb{B}^+_\dR)}(\ast_{\mu'})=\Gr_{[\mu']}\]
    fitting in a commutative diagram 

\begin{equation}\label{eq.comm-diag-orbits}\begin{tikzcd}
	{G(\mathbb{B}^+_\dR)/\mathrm{Stab}_{G(\mathbb{B}^+_\dR)}(\ast_{\mu})} & {\GL_n(\mathbb{B}^+_\dR)/\mathrm{Stab}_{\GL_n(\mathbb{B}^+_\dR)}(\ast_{\mu'})} \\
	{\Gr_{[\mu]}} & {\Gr_{[\mu']}}
	\arrow[from=1-1, to=1-2]
	\arrow[from=1-1, to=2-1]
	\arrow["\sim"{description}, from=1-2, to=2-2]
	\arrow[hook, from=2-1, to=2-2]
\end{tikzcd}\end{equation}
By \cref{prop.closed-schubert-functor}, the map $\Gr_{[\mu]}\hookrightarrow \Gr_{[\mu']}$ induced by $\rho$ that appears as the bottom horizontal arrow in \cref{eq.comm-diag-orbits} is a closed immersion. In particular, since both $G(\mathbb{B}^+_\dR)/\mathrm{Stab}_{G(\mathbb{B}^+_\dR)}(\ast_\mu)$ and $\Gr_{[\mu]}$ have the same geometric points, it suffices to check that the top horizontal arrow is also a closed immersion. Indeed, this will imply the left vertical arrow is quasi-compact and then we may apply \cite[Lemma 12.5]{Scholze.EtaleCohomologyOfDiamonds} to conclude that it is an isomorphism. 

It thus remains to show that the natural map
\begin{equation}\label{eq.natural-inclusion-map} G(\mathbb{B}^+_\dR)/\mathrm{Stab}_{G(\mathbb{B}^+_\dR)}(\ast_{\mu}) \rightarrow \GL_n(\mathbb{B}^+_\dR)/\mathrm{Stab}_{\GL_n(\mathbb{B}^+_\dR)}(\ast_{\mu'}) \end{equation}
is a closed immersion. To that end, we write $G'=\GL_n$, $H=\mathrm{Stab}_{G(\mathbb{B}^+_\dR)}(\ast_{\mu})$, and $H'=\mathrm{Stab}_{G'(\mathbb{B}^+_\dR)}(\ast_{\mu'})$, so that $H=H' \cap G$. 

We use the notation for the filtrations by principal congruence subgroups as in \cref{lemma.stabilizer-subquotient-structure}. The map appearing in \cref{eq.natural-inclusion-map} can be written in this notation as $G_0/H_0 \rightarrow G'_0/H'_0$ and is filtered by the maps
\[ G_{i}\cdot H_0 / H_0 \rightarrow G'_i\cdot H'_0 / H'_0 \textrm{ for $i \geq 0$.} \]
The subquotient maps for this filtration are maps 
\begin{equation}\label{eq.subqmap} G_i\cdot H_0/G_{i+1}\cdot H_0 \rightarrow G'_i\cdot H'_0/G'_{i+1}\cdot H'_0.\end{equation}
We show $G_0/H_0 \rightarrow G'_0/H'_0$ is a closed immersion by showing that each of the subquotient maps \cref{eq.subqmap} is a closed immersion. We rewrite
\begin{multline*} G_i \cdot H_0 / G_{i+1} \cdot H_0 = G_i/ \left((G_{i+1} \cdot H_0) \cap G_i\right) = (G_i/G_{i+1})/ \left(G_{i+1} \cdot (H_0 \cap G_i) / G_{i+1} \right)\\
 = (G_i/G_{i+1})/ \left( H_0 \cap G_i / H_0 \cap G_{i+1} \right) = (G_i/G_{i+1}) /  (H_i/H_{i+1}) \end{multline*}
and similarly for $G_i'$. For $i=0$, using \cref{lemma.stabilizer-subquotient-structure}, \cref{eq.subqmap} thus becomes
\[ G(\mathcal{O})/P_{\mu^{-1}}(\mathcal{O}) \rightarrow G'(\mathcal{O})/P_{{\mu'}^{-1}}(\mathcal{O}), \]
which is the closed immersion of associated flag varieties $\Fl_{\mu^{-1}}^\diamond \rightarrow \Fl_{{\mu'}^{-1}}^\diamond.$ Similarly, for $i >0$, \cref{eq.subqmap} becomes the inclusion of a vector subspace
\[ (\Lie G)\{i\}/ \Fil_{\mu^{-1}}^{-i}(\Lie G)\{i\} \rightarrow (\Lie G')\{i\}/ \Fil_{{\mu'}^{-1}}^{-i}(\Lie G')\{i\}, \]
and thus is a closed immersion. 
\end{proof}

\subsection{The general case}\label{ss.general-case}

The proof of \cref{prop.red-case} does not generalize to non-reductive connected linear algebraic groups since it uses that, for a reductive group $G$, any embedding $G \rightarrow \GL_n$ realizes $\Gr_G$ as a closed subsheaf of $\Gr_{\GL_n}$. For a non-reductive group, one can at best get a locally closed embedding $\Gr_G \rightarrow \Gr_{\GL_n}$ by choosing $G \rightarrow \GL_n$ realizing $G$ as an observable subgroup (see \cite[Section 5.1]{HoweKlevdal.AdmissiblePairsAndpAdicHodgeStructuresIITheBiAnalyticAxLindemannTheorem}). As a locally closed embedding may not be quasi-compact, this breaks the step where one applies \cite[Lemma 12.5]{Scholze.EtaleCohomologyOfDiamonds}. Instead, we reduce the general case to the reductive case by making arguments similar to those appearing in \cite[Section 5.2]{HoweKlevdal.AdmissiblePairsAndpAdicHodgeStructuresIITheBiAnalyticAxLindemannTheorem}.

\begin{proof}[Proof of \cref{theorem.action-transitive} in the general case]
We fix a Levi decomposition $G=MU$. By base change to a finite extension $L'/L$, we may assume $M$ is split so that $L=L([\mu])$ and we can moreover fix a representative $\mu:\mathbb{G}_{m,L} \rightarrow G$ that factors through a cocharacter $\mu_M$ of $M$. 

By push-out from $M$ to $G$ and from $G$ to $M$, we obtain maps \[ \Gr_M \rightarrow \Gr_G \rightarrow \Gr_M \]
that compose to the identity and are equivariant along the natural maps $M(\mbb{B}_\dR) \rightarrow G(\mbb{B}_\dR) \rightarrow M(\mbb{B}_\dR)$. These restrict to maps 
\begin{equation}\label{eq.v-sheaf-schubert-cell-retraction} \Gr_{[\mu_M]} \rightarrow \Gr_{[\mu]} \rightarrow \Gr_{[\mu_M]} \end{equation}
that are equivariant along the natural maps $M(\mbb{B}^+_\dR) \rightarrow G(\mbb{B}^+_\dR) \rightarrow M(\mbb{B}^+_\dR)$, and moreover $\ast_{\mu_M} \mapsto \ast_{\mu} \mapsto \ast_{\mu_M}.$ 

It follows from the reductive case treated already above that $m \mapsto m \cdot \ast_{\mu_M}$ is a surjection  $M(\mbb{B}^+_\dR) \rightarrow \Gr_{[\mu_M]} $ in the \'{e}tale topology. To establish the general result, we need to show that $g \mapsto g \cdot \ast_{\mu}$ is a surjection $G(\mbb{B}^+_\dR) \rightarrow \Gr_{[\mu]}$ in the \'{e}tale topology.

So, suppose $x \in \Gr_{[\mu]}(\Spa(R,R^+)) \subseteq \Gr_G(\Spa(R,R^+))$. After passing to an \'{e}tale cover, by the definition of $\Gr_G$, we may assume $x=g \cdot \ast_1$ for $g \in G(\mbb{B}_\dR(R))$. We write $g=u m$ for $u \in U(\mbb{B}_\dR(R))$ and $m \in M(\mbb{B}_\dR(R))$. The image of $x = g \cdot \ast_1= um \cdot \ast_1$ in $\Gr_{M}(\Spa(R,R^+))$ is equal to $m \cdot \ast_{1_M}$ (because $u$ maps to the identity element in $M$) and, since $x \in \Gr_{[\mu]}(\Spa(R,R^+))$, we deduce that $m \cdot \ast_{1_M}$ is in $\Gr_{[\mu_M]}(\Spa(R,R^+))$. Thus, by the surjectivity for $M$, after passing to another \'{e}tale cover we may assume $m = m_1 \xi^\mu m_2$ for $m_1,m_2 \in M(\mbb{B}^+_\dR)$. It follows that 
\[ m_1^{-1}g= m_1^{-1} u m_1 \xi^\mu m_2 = u' \xi^\mu m_2\]
for $u'=m_1^{-1} u m_1 \in U(\mbb{B}_\dR)$, and thus 
\[ m_1^{-1}x=m_1^{-1}g \cdot \ast_1 = u' \xi^\mu m_2 \cdot \ast_1 = u' \xi^\mu \cdot \ast_1 =u' \cdot \ast_{\mu}. \]
Thus $m_1 u' \cdot \ast_{\mu}=x$, so it suffices to show that we can replace $u'$ here with an element of $U(\mbb{B}^+_\dR)$. To accomplish this we will use the fact that $u' \cdot \ast_{\mu} = m_1^{-1} x \in \Gr_{\mu}(\Spa(R,R^+))$. As in \cite[Lemma 3.4.3 and preceding paragraph]{HoweKlevdal.AdmissiblePairsAndpAdicHodgeStructuresITranscendenceOfTheDeRhamLattice} (note that our $U$ becomes $U'$ in loc. cit., and is contained in the unipotent subgroup of a Borel in $\GL_n$, denoted $U$ in loc. cit.), we may consider the subgroups $U_{>0},$ $U_{<0},$ and $U_{0}$ of $U$ whose Lie algebras consist of those character spaces on which $\mu$ acts positive, negatively, or by zero, and multiplication then induces an isomorphism $U=U_{>0}\times U_0 \times U_{<0}$. Writing $u'=u'_{>0} u'_0 u'_{<0}$ according to this decomposition, we claim it follows from \cite[Lemma 3.4.3]{HoweKlevdal.AdmissiblePairsAndpAdicHodgeStructuresITranscendenceOfTheDeRhamLattice} that $u'_{>0} \cdot \ast_{[\mu]}=u' \cdot \ast_{[\mu]}$ and that $u'_{>0} \in U_{>0}(\mbb{B}^+_\dR)$. Explicitly, by fixing an embedding of $G$ into $\GL_n$, we can apply  \cite[Lemma 3.4.3]{HoweKlevdal.AdmissiblePairsAndpAdicHodgeStructuresITranscendenceOfTheDeRhamLattice} (noting again that our $U$ becomes $U'$ in loc. cit., and is contained in the unipotent subgroup of a Borel in $\GL_n$, denoted $U$ in loc. cit.) to see that $u'_{>0} \in U(\mbb{B}^+_\dR)$, and $\xi^{-\mu}u'_0 u'_{<0}\xi^\mu \in  U(\mbb{B}^+_\dR) $. Precisely, \cite[Lemma 3.4.3]{HoweKlevdal.AdmissiblePairsAndpAdicHodgeStructuresITranscendenceOfTheDeRhamLattice} gives this at each geometric point, but each of these is a closed condition on $\Spa(R,R^+)$ by \cite[Lemma 2.5]{HoweKlevdal.AdmissiblePairsAndpAdicHodgeStructuresIITheBiAnalyticAxLindemannTheorem}. It then follows that 
\[ u' \cdot \ast_{\mu}=u' \xi^\mu \cdot \ast_1 = u'_{>0} \xi^{\mu} (\xi^{-\mu} u'_{0}u'_{<0} \xi^{\mu}) \cdot \ast_1 = u'_{>0} \xi^{\mu} \cdot \ast_1 = u'_{>0} \cdot \ast_{\mu}. \]
\end{proof}

\bibliographystyle{plain}
\bibliography{references, preprints}

@article{CesnaviciusYoucis.TheAnalyticTopologySufficesForTheBplusDrGrassmannian,
      title={The analytic topology suffices for the {$B_{\mathrm{dR}}^+$}-grassmannian}, 
      author={Kestutis Cesnavicius and Alex Youcis},
      year={2024},
      journal={To appear in the Proceedings of the Simons Symposium on p-adic Hodge theory, ar{X}iv:2303.11710},
      primaryClass={math.AG},
      url={https://arxiv.org/abs/2303.11710}, 
}

@article{Howe.InscriptionTwistorsAndPAdicPeriods,
AUTHOR={Howe, Sean},
TITLE={Inscription, twistors, and $p$-adic periods},
YEAR={2025},
JOURNAL={ar{X}iv:2508.11589}
}

@misc{stacks-project,
    shorthand    = {Stacks},
    author       = {The {Stacks Project Authors}},
    title        = {\textit{Stacks Project}},
    howpublished = {\url{https://stacks.math.columbia.edu}},
    year         = {2018},
  }

@article {Kottwitz.ShimuraVarietiesAndTwistedOrbitalIntegrals,
    AUTHOR = {Kottwitz, Robert E.},
     TITLE = {Shimura varieties and twisted orbital integrals},
   JOURNAL = {Math. Ann.},
  FJOURNAL = {Mathematische Annalen},
    VOLUME = {269},
      YEAR = {1984},
    NUMBER = {3},
     PAGES = {287--300},
      ISSN = {0025-5831,1432-1807},
   MRCLASS = {11G18 (11F70 11G40 22E50)},
  MRNUMBER = {761308},
MRREVIEWER = {Ernst-Wilhelm\ Zink},
       DOI = {10.1007/BF01450697},
       URL = {https://doi.org/10.1007/BF01450697},
}

@Article{KedlayaLiu.RelativepAdicHodgeTheoryFoundations,
  author     = {Kedlaya, Kiran S. and Liu, Ruochuan},
  journal    = {Ast\'{e}risque},
  title      = {Relative {$p$}-adic {H}odge theory: foundations},
  year       = {2015},
  issn       = {0303-1179},
  number     = {371},
  pages      = {239},
  fjournal   = {Ast\'{e}risque},
  isbn       = {978-2-85629-807-7},
  mrclass    = {14C30 (14G22)},
  mrnumber   = {3379653},
  mrreviewer = {Giovanni Rosso},
}

@article {Rapoport.APositivityPropertyOfTheSatakeIsomorphism,
    AUTHOR = {Rapoport, Michael},
     TITLE = {A positivity property of the {S}atake isomorphism},
   JOURNAL = {Manuscripta Math.},
  FJOURNAL = {Manuscripta Mathematica},
    VOLUME = {101},
      YEAR = {2000},
    NUMBER = {2},
     PAGES = {153--166},
      ISSN = {0025-2611,1432-1785},
   MRCLASS = {22E50 (14F30 20G25)},
  MRNUMBER = {1742251},
MRREVIEWER = {Adolfo\ Quir\'os},
       DOI = {10.1007/s002290050010},
       URL = {https://doi.org/10.1007/s002290050010},
}

@Article{Scholze.pAdicHodgeTheoryForRigidAnalyticVarieties,
  author     = {Scholze, Peter},
  journal    = {Forum Math. Pi},
  title      = {{$p$}-adic {H}odge theory for rigid-analytic varieties},
  year       = {2013},
  issn       = {2050-5086},
  pages      = {e1, 77},
  volume     = {1},
  doi        = {10.1017/fmp.2013.1},
  fjournal   = {Forum of Mathematics. Pi},
  mrclass    = {14G22 (14C30 14F30 14G20 32J27 32P05)},
  mrnumber   = {3090230},
  mrreviewer = {Hui June Zhu},
  url        = {https://doi-org.stanford.idm.oclc.org/10.1017/fmp.2013.1},
}

@Book{ScholzeWeinstein.BerkeleyLecturesOnPAdicGeometryAMS207,
  author    = {Peter Scholze and Jared Weinstein},
  publisher = {Princeton University Press},
  title     = {Berkeley Lectures on p-adic Geometry: (AMS-207)},
  year      = {2020},
  isbn      = {9780691202082},
  url       = {http://www.jstor.org/stable/j.ctvs32rc9},
}

@Article{FarguesScholze.GeometrizationOfTheLocalLanglandsCorrespondence,
  author  = {Fargues, Laurent and Scholze, Peter},
  journal = {To appear in Asterisque, ar{X}iv:2102.13459},
  title   = {Geometrization of the local {L}anglands correspondence},
}

@Article{HoweKlevdal.AdmissiblePairsAndpAdicHodgeStructuresITranscendenceOfTheDeRhamLattice,
  author = {Howe, Sean and Klevdal, Christian},
  title  = {Admissible pairs and $p$-adic {H}odge structures {I}: {T}ranscendence of the de {R}ham lattice},
  journal={To appear in Algebra and Number Theory, ar{X}iv:2308.11065},
year={2023}
}

@Article{HoweKlevdal.AdmissiblePairsAndpAdicHodgeStructuresIITheBiAnalyticAxLindemannTheorem,
  author = {Howe, Sean and Klevdal, Christian},
  title  = {Admissible pairs and $p$-adic {H}odge structures {II}: The bi-analytic {A}x-{L}indemann theorem},
  journal={Inventiones Mathematicae},
year={2025}
}

@Article{Scholze.EtaleCohomologyOfDiamonds,
  author    = {Scholze, Peter},
  title     = {Etale cohomology of diamonds},
  year      = {2017},
  journal = {To appear in Asterisque, arXiv:1709.07343}
}
\end{document}